\documentclass[letterpaper,12pt,leqno]{article}

\usepackage{fullpage}

\usepackage{hyperref}

\usepackage{amsmath}
\usepackage{amsfonts}
\usepackage{amssymb}
\usepackage{tabularx}
\usepackage{color}
\usepackage{fancyhdr}
\usepackage{tikz-cd}
\usepackage{enumitem}

\def\theoremparent{section}
\usepackage{notes}

\numberwithin{equation}{theorem}

\def\A{\mathbb{A}}

\def\cH{\mathcal{H}}

\def\d{\partial}

\def\Der{\mathrm{Der}}

\def\EE{\mathbb E}

\def\g{\mathfrak{g}}

\def\h{\mathfrak{h}}
\def\rmht{\mathrm{ht}}

\def\k{\mathbf{k}}

\def\proof{\noindent{\em Proof:}\ }

\def\qed{\hfill\lower 1em\hbox{$\square$}\vskip 1em}

\def\SS{\mathbb{S}}

\def\sk{\mathop{\mathrm{sk}}\nolimits}

\def\toto{\text{\ \raise0.2em\hbox to 0pt {$\to$}\lower0.2em\hbox{$\to$}\ }}

\def\V{\mathcal{V}}

\begin{document}

\title{Gelfand-Fuchs cohomology for affine superspaces $\A^{m,n}$}
\author{Slava Pimenov}
\date{\today}
\titlepage
\maketitle

\tableofcontents
\vskip 5em

\setcounter{section}{-1}
\section{Introduction}
Let $\A^{m,n}$ be the the affine super space of even dimension $m$ and odd dimension $n$ over an algebraically closed field $\k$ of characteristic $0$.
Consider Lie superalgebras $\V_{m,n}$ of vector fields in the formal neighborhood of $0 \in \A^{m,n}$.  This is a topological Lie superalgebra with $\k$-linear
topology induced by the defining ideal of point $0 \in \A^{m,n}$. We are interested in the continuous cohomology groups $H^\bullet(\V_{m,n}, \k)$.

Previously established results cover cases $0 \le m \le n$, as well as $n = 0$ and $n = 1$ and arbitrary $m \ge 0$ (\cite{Gel-Fuchs}, \cite{Fuks}, \cite{Koszul}, \cite{Fuks-super}, \cite{gfn1}).
These results are collected below in theorem \ref{thm_known_gf}. The main result of this paper is the following theorem that has been stated as a conjecture in the previous paper (\cite{gfn1}).

\begin{atheorem}
\label{thm_intro_gfmn}
For any $m \ge n \ge 0$ we have an isomorphism
$$
H^\bullet(\V_{m,n}, \k) \ \isom\ H^\bullet(\SS^{2n} X_{2(m - n)}, \k).
$$
\end{atheorem}
Here $\SS$ denotes the topological suspension functor, and $X_{2(m-n)}$ is the pullback of the tautological $GL(m, \CC)$-torsor over $BGL(m, \CC)$ to the
$2(m-n)$-dimensional skeleton of $BGL(m, \CC)$ consisting of cells of dimensions up to $2(m - n)$.

Combined with the previously established results this completely settles the question of local Gelfand-Fuchs cohomology for super-manifolds.

\begin{nparagraph}
The main tool in the calculation is the following theorem regarding cohomology of Lie superalgebras $\gl(m, n)$. Let $V$ be the standard representation of $\gl(m, n)$, and denote
$\Sigma^\lambda(V)$ the Schur functor corresponding to a diagram $\lambda$. We will write $\cH_{m,n}$ for the set of diagrams contained in a thick hook with $m$ rows and $n$ columns.

\begin{atheorem}
\label{thm_intro_glmn}
Let $\g = \gl(m, n)$ with $m \ge n \ge 0$, and
$$
\lambda \in \cH_{m - n + k, k} - \cH_{m - n + k - 1, k - 1}
$$
for some $0 \le k \le n$. Then
$$
H^\bullet(\g, \Sigma^\lambda(V) \tensor \Sigma^\lambda(V^*)) \ \isom\ \k[e_1, \ldots, e_{2m - 1}] \tensor \k[e'_{2(n - k) + 1}, \ldots, e'_{2n - 1}].
$$
\end{atheorem}
This theorem appears to be a new result and may be of interest beyond the Gelfand-Fuchs cohomology theory.

\end{nparagraph}

\begin{nparagraph}[Outline of the paper.]
In section \ref{sec_notations} we recall the relevant notations and results that are used in this paper.

The section \ref{sec_glmn} is dedicated to the proof of theorem \ref{thm_intro_glmn}. We proceed by induction on the number of odd variables $n$ and use the spectral
sequence for the Lie subalgebra
$$
\gl(m, n - 1) \oplus \gl(1) \into \gl(m, n)
$$
to reduce the question to $\gl(m, n - 1)$.

First we observe that the first layer of this spectral sequence has a universal structure, which allows us to compare spectral sequences for a fixed diagram $\lambda$
but different values of $m$ and $n$. We combine this with the special case of theorem \ref{thm_intro_glmn} for $n = 1$ that was established in (\cite{gfn1}) to identify
all the diagrams contributing to the second layer of the spectral sequence. Then by direct examination of the second and third layers we establish the required isomorphism.

We would like to point out that this does not provide an independent proof of theorem \ref{thm_intro_glmn} for $n = 1$ as we use this result to greatly simplify analysis
of the first layer of the spectral sequence. Combining the proof in this paper and the proof of the special case for $n = 1$, the overall process
looks as follows. We start with $\gl(1, 1)$ as the base of induction, then we grow number of even variables to get to $\gl(m, 1)$ then we grow number of odd variables
and arrive to $\gl(m, n)$.

In the case of $\gl(1, 1)$ we have a description of the indecomposable components of the coefficient module $\Sigma^\lambda(V) \tensor \Sigma^\lambda(V^*)$
that give rise to the cohomology groups in theorem \ref{thm_intro_glmn}. Presently, we do not have a similar description for other $\gl(n, n)$.

\end{nparagraph}

\begin{nparagraph}
The section \ref{sec_Vmn} deals with the proof of theorem \ref{thm_intro_gfmn}. It follows the general process developed in (\cite{gfn1}) which in turn is a refinement of the proof of Gelfand and Fuchs in the classical case.
We consider the spectral sequence for the Lie subalgebra $\gl(m, n) \into \V_{m,n}$ and identify all the diagrams contributing to its first layer. We observe that
up to transposition these diagrams are the same as in the spectral sequence for $\gl(n - 1, m + 1) \into \V_{n - 1, m + 1}$, which allows us to compare these two
spectral sequences.

In the latter case we have $n - 1 < m + 1$, so this is covered by the previous result of Astashkevich and Fuchs (\cite{Fuks-super}) which says that the cohomology
$H^\bullet(\V_{n - 1, m + 1}, \k)$ is isomorphic to the cohomology of $(2m + 1)$-dimensional sphere. This makes the analysis of the original spectral sequence for
$\V_{m,n}$ much simpler and allows us to establish the required isomorphism.
\end{nparagraph}

\vskip 1em
The author would like to thank BIMSA (Beijing Institute of Mathematical Sciences and Applications) for providing excellent working conditions during preparation of this paper.

\vskip 5em
\section{Notations and recollections}
\label{sec_notations}
We will retain conventions and notations from \cite{gfn1}. Here we will briefly recall them and state the relevant results that will be used in this
paper.

\begin{nparagraph}[Young diagrams and Schur functors.]
Let $\lambda$ be a Young diagram of size $d$, in other words it is an unordered partition $\lambda = (\lambda_1, \ldots, \lambda_k)$ of $d$, with
$\lambda_1 \ge \lambda_2 \ge \ldots \ge \lambda_k > 0$ and $\sum_{i=1}^k \lambda_i = d$. We will refer to $k$ as the height of $\lambda$ and write
$k = \rmht(\lambda)$, and $d = |\lambda|$. We will denote $\lambda'$ the transposed Young diagram, specifically we put $\lambda'_j = \max \{ i \mid \lambda_i \ge j\}$.

For any diagram $\lambda$ we construct a truncated diagram $\wbar\lambda$ obtained from $\lambda$ by removing the first column. In other words we put
$\wbar\lambda_i = \max \{\lambda_i - 1, 0\}$. Furthermore, we construct an extended diagram $\wtilde\lambda$ by adding to $\lambda$ the first column
of height $d = |\lambda|$. Formally, we put $\wtilde\lambda_i = \lambda_i + 1$ for $1 \le i \le d$.

For any $m, n \ge 0$ we consider a subset $\cH_{m,n}$ of Young diagrams of arbitrary size, consisting of diagrams contained in a thick hook with $m$ rows and $n$ columns.
More precisely $\lambda \in \cH_{m,n}$ if and only if $\lambda_i \le n$ whenever $i > m$. By convention, if either $m$ or $n$ is negative we put $\cH_{m,n}$
to be an empty set.

For any partition $\lambda$ we will denote $\Sigma^\lambda$ the corresponding Schur functor acting on the symmetric monoidal category of super vector spaces.
For a super vector space $V = (V_0, V_1)$ of dimension $(m, n)$, with $m = \dim V_0$ and $n = \dim V_1$, the Schur functor $\Sigma^\lambda(V)$ is non-zero
if and only if $\lambda \in \cH_{m,n}$.

We will denote by $S^n$ and $\Lambda^n$ the functors of symmetric and exterior powers respectively and recall that for any two super vector spaces $V$ and $W$
we have isomorphisms of $(\gl(V) \times \gl(W))$-modules
\begin{equation}
\label{equ_sym_ext}
S^n(V \tensor W) = \bigoplus_{|\lambda| = n} \Sigma^\lambda(V) \tensor \Sigma^\lambda(W),\quad\quad
\Lambda^n(V \tensor W) = \bigoplus_{|\lambda| = n} \Sigma^\lambda(V) \tensor \Sigma^{\lambda'}(W).
\end{equation}

\end{nparagraph}

\begin{nparagraph}[Lie superalgebra $\gl(m,n)$.]
Let $V$ be a super vector space of dimension $(m, n)$, then we write $\gl(m,n)$ for the Lie superalgebra of endomorphisms $\End(V) \isom V \tensor V^*$.
We will refer to $V$ as the standard representation of $\gl(m,n)$. We would like to point out here that even though Lie superalgebras $\gl(m,n)$ and
$\gl(n,m)$ are isomorphic, their standard representations are different. Specifically, the standard representation $W$ of $\gl(n,m)$ is obtained
from $V$ by the change of parity $W \isom \Pi(V)$. Therefore, up to change of parity we have isomorphisms $\Sigma^\lambda(W) = \Sigma^\lambda(\Pi V) \isom \Sigma^{\lambda'}(V)$.

Let $\g = \gl(m,n)$, and $V$ its standard representation, we are interested in the cohomology spaces $H^\bullet(\g, \Sigma^\lambda(V) \tensor \Sigma^\lambda(V^*))$.
In \cite{gfn1} we have established the following result.

\end{nparagraph}

\begin{theorem}
\label{thm_cohom_glm1}
Let $\g = \gl(m,1)$, $V$ the standard representation of $\g$ and $\lambda \in \cH_{m,1}$, then
$$
H^\bullet(\g, \Sigma^\lambda(V) \tensor \Sigma^\lambda(V^*)) = \begin{cases}
\k[e_1, \ldots, e_{2m-1}],&\text{if $\rmht(\lambda) \le m - 1$,}\\
\k[e_1, \ldots, e_{2m-1}, e'_1],&\text{otherwise},
\end{cases}
$$
where generators $e_i$ are of cohomological degree $i$ and $e'_1$ is of degree $1$.
\end{theorem}

Notice that the condition $\rmht(\lambda) \le m - 1$ here can be rewritten as $\lambda \in \cH_{m - 1,0}$.

\begin{nparagraph}[Lie superalgebra $\V_{m,n}$.]
Consider a (super)commutative superalgebra
$$
\O_{\A^{m,n}} = \k[x_1, \ldots x_m, \xi_1, \ldots \xi_n]
$$
of algebraic functions on the affine superspace $\A^{m,n}$.
We assume that variables $x_i$ are even and $\xi_j$ are odd. We denote $\hat\O_{\A^{m,n}} = \k[[x_1, \ldots, x_m, \xi_1, \ldots, \xi_n]]$ its completion at zero,
equipped with the inverse limit topology. We are interested in the Lie superalgebra of continuous derivations
$$
\V_{m,n} = \Der_{\mathrm{cont}} (\hat\O_{\A^{m,n}}).
$$

Explicitly, $\V_{m,n}$ is formed by elements $\sum f_i {\d_{x_i}} + \sum g_j {\d_{\xi_j}}$, with $f_i, g_j \in \hat\O_{\A^{m,n}}$. The bracket is given by
the action of derivations $\d_{x_i}$ and $\d_{\xi_j}$ on functions. It contains $\gl(m,n)$ as a subalgebra spanned by elements with linear coefficients:
$\{x_i \d_{x_j}, x_i \d_{\xi_j}, \xi_i \d_{x_j}, \xi_i \d_{\xi_j}\}$.

We will write $H^\bullet(\V_{m,n}, \k)$ for the continuous cohomology spaces with respect to topology on $\V_{m,n}$ induced by the inverse limit topology on $\hat\O_{\A^{m,n}}$.
We recall the previously established results regarding this cohomology, for details we refer to \cite{Fuks}, \cite{Fuks-super}, \cite{gfn1}.

The cohomology of $\V_{m,n}$ will be related to the cohomology of various topological spaces, that can be constructed using the following procedure.
Consider the topological group $GL(m, \CC)$, and let
$BGL(m)$ be its classifying space. Denote by $p\from EGL(m) \to BGL(m)$ the tautological principal $GL(m)$-bundle over the classifying
space. Let us write $\sk_d BGL(m)$ for the $d$-dimensional skeleton of $BGL(m)$, i.e. the subspace formed by all cells of
dimension up to $d$. The spaces that will be of interest to us are $X_d = p^{-1}(\sk_d BGL(m))$ for various $d$, i.e. the pullbacks of the tautological bundle to the $d$-dimensional skeleta.
Furthermore, denote by $\SS$ the topological suspension functor.

\end{nparagraph}

\begin{theorem}
\label{thm_known_gf}
We have the following isomorphisms.
\begin{enumerate}[label=\alph*)]
\item{\cite{Fuks}} If $n = 0$, then
$$
H^\bullet(\V_{m,0}, \k) \ \isom\  H^\bullet(X_{2m}, \k).
$$

\item{\cite{Fuks-super}} If $m < n$, then
$$
H^\bullet(\V_{m,n}, \k) \isom H^\bullet(S^{2n-1}, \k).
$$

\item{\cite{Fuks-super}} If $m = n$, then
$$
H^\bullet(\V_{n,n}, \k) \ \isom\ H^\bullet(\SS^{2n} GL(n, \CC), \k).
$$

\item{\cite{gfn1}} If $n = 1$, then
$$
H^\bullet(\V_{m,1}, \k) \ \isom\ H^\bullet(\SS^2 X_{2(m-1)}, \k).
$$
\end{enumerate}
\end{theorem}

\begin{nparagraph}[Spectral sequence.]
\label{par_sp_seq}
The main tool used in calculations in this paper is the spectral sequence relating the cohomology of a Lie superalgebra $\g$ and its subalgebra $\h$. For any
$\g$-module $M$ we have an increasing filtration of the chain complex $\Lambda^\bullet \g \tensor M$ by the number of elements from $\h$. It induces a
decreasing filtration on the cochain complex $C^\bullet(\g, M)$, giving rise to a spectral sequence
$$
E_1^{pq} = H^q(\h, \ \Hom(\Lambda^p (\g / \h), M)\ ) \Rightarrow H^{p+q}(\g, M).
$$

We will use cohomological indexing convention for the spectral sequence: on the layer $E_r$ we have differentials
$$
d_r\from E_r^{p,q} \to E_r^{p+r, q-r+1}.
$$

\end{nparagraph}

\vskip 5em
\section{Cohomology of $\gl(m,n)$}
\label{sec_glmn}
Let $\g = \gl(m,n)$, and $V$ its standard representation. Dimension of the super vector space $V$ is $(m, n)$. Consider a direct sum decomposition of $V$
into two subspaces $V = W \oplus E$, such that $\dim W = (m, n - 1)$ and $\dim E = (0, 1)$. Denote by $\h$ the subalgebra of $\g$ that preserves this
decomposition, in other words
$$
\h = \End(W) \oplus \End(E) \isom \gl(m, n - 1) \oplus \gl(1) \into \g.
$$
We identify $W$ with the standard representation of $\gl(m, n - 1)$ and $E$ with the standard representation of $\gl(0, 1) \isom \gl(1)$. The quotient
space $\g / \h$ is isomorphic to $W \tensor E^* \oplus W^* \tensor E$. We have

$$
\Lambda^p(\g / \h) \ \isom\ \bigoplus_{i + j = p} \Lambda^i(W \tensor E^*) \tensor \Lambda^j(W^* \tensor E).
$$
Since $E$ is of odd dimension $1$, we can rewrite this as follows

\begin{equation}
\label{equ_VEW_ext}
\Lambda^p(\g / \h) \ \isom\ \bigoplus_{i + j = p} S^i(W) \tensor S^j(W^*) \tensor \Lambda^i(E^*) \tensor \Lambda^j(E).
\end{equation}

\begin{nparagraph}
We are interested in the cohomology of $\g$ with coefficients in $\Sigma^\lambda(V) \tensor \Sigma^\lambda(V^*)$. Using decomposition $V = W \oplus E$
we can write

\begin{equation}
\label{equ_VEW_coeff}
\Sigma^\lambda(V) \ \isom\ \bigoplus_{\mu} \Sigma^\mu(W) \tensor \Lambda^p(E),
\end{equation}
where the sum is taken over all diagrams $\mu$ obtained from $\lambda$ by removing at most one box from each row, and $p = |\lambda| - |\mu|$ is the total number
of removed boxes. We also have a similar expansion for $\Sigma^\lambda(V^*)$.

Let us consider the spectral sequence $E$ for the Lie subalgebra $\h \into \g$ and coefficients $\Sigma^\lambda(V) \tensor \Sigma^\lambda(V^*)$:
$$
E_1^{pq} = H^\bullet(\gl(m, n - 1) \oplus \gl(1), \Lambda^p(\g / \h)^* \tensor \Sigma^\lambda(V) \tensor \Sigma^\lambda(V^*)).
$$
Combining \ref{equ_VEW_ext} and \ref{equ_VEW_coeff} we find that diagrams $\mu$ contributing to the first layer of this spectral sequence are obtained
from $\lambda$ by first removing say $i$ boxes in such a way that from each row we remove at most one box, and then adding say $j$ boxes in such a way
that in each column we add at most one box. In particular, we immediately see that the diagrams $\mu$ appearing in $E_1$ are of height at most
$\rmht(\lambda) + 1$.

Furthermore, the weight with respect to the action of the subalgebra $\gl(1) \into \h$ of a component corresponding to a diagram $\mu$ is $|\lambda| - |\mu|$, in other words
it depends only on the diagram $\mu$ itself, and not on the specific way it was obtained from $\lambda$ by the procedure described above. We have a similar picture
on the dual side with components containing $\Sigma^\nu(W^*)$, however
since the cohomology $H^\bullet(\gl(m, n - 1), \Sigma^\mu(W) \tensor \Sigma^\nu(W^*))$ vanishes unless $\mu = \nu$ it is sufficient to keep track only
of the diagrams $\mu$. We will refer to such components as components of type $\mu$.

\end{nparagraph}

\begin{nparagraph}
The differential on the first layer of the spectral sequence is induced by maps between the coefficients of $H^\bullet(\h, -)$, which in turn corresponds to the action of
$\g / \h$ on the coefficients $\Sigma^\lambda(V) \tensor \Sigma^\lambda(V^*)$. In terms of the above decomposition
a differential between two components of type $\mu$ corresponds to two different ways of obtaining diagram $\mu$ from $\lambda$. More precisely, let $\nu$ be the
intermediate diagram in the process of constructing $\mu$ from $\lambda$, i.e. $\nu$ is a common subdiagram of $\mu$ and $\lambda$ such that
$\lambda - \nu$ has at most one box in each row, and $\mu - \nu$ has at most one box in each column.

We say that a box of the diagram $\nu$ is {\em flippable} if the diagram $\nu_1$ obtained from $\nu$ by removing this box provides another valid way of obtaining
$\mu$ from $\lambda$. The differential in $E_1$ is a linear combination of maps into components obtained by flipping a box (either for $\Sigma^\mu(W)$ or for $\Sigma^\mu(W^*)$).
It is clear from the construction that for each $\mu$ the subcomplex of components of type $\mu$ is bounded.

\end{nparagraph}

\begin{nparagraph}[Universal complex.]
The previous discussion of the first layer of the spectral sequence can by summarized by saying that it is ``universal'' in a certain sense. We will make this statement
more precise. Let $S_d$ be a symmetric group on $d$ elements and $S_\bullet$ denote the collection of all $S_d$ for $d \ge 0$. An $S_\bullet$-module $M$ is a direct sum
$M = \bigoplus_{d \ge 0} M_d$, where each $M_d$ is an $S_d$-module. Since the category of $S_\bullet$-modules is semisimple, and simple modules correspond to partitions $\lambda$
of arbitrary size, we can further decompose
$$
M \ =\ \bigoplus_\lambda L_\lambda \tensor \Hom_{S_\bullet}(L_\lambda, M),
$$
where $L_\lambda$ is the simple module corresponding to a partition $\lambda$. To simplify notation we will write $M_\lambda = \Hom_{S_\bullet}(L_\lambda, M)$.

For a super vector space $V$ the Schur functor $\Sigma(M, V)$ is defined by
$$
\Sigma(M, V) \ =\ \bigoplus_{d \ge 0} M_d \tensor_{\k[S_d]} V^{\tensor d} \ \isom\ \bigoplus_\lambda \Sigma^\lambda(V) \tensor M_\lambda,
$$
where $S_d$ acts on the tensor power $V^{\tensor d}$ by permuting factors. Of course, as was already mentioned in the previous section, the dimension of $V$ imposes restriction
on what terms contribute to the direct sum. If $\dim V = (m, n)$, then the contribution comes only from diagrams $\lambda \in \cH_{m,n}$.

The above discussion of the first layer of the spectral sequence can be stated as the following lemma.

\end{nparagraph}

\begin{lemma}
\label{lemma_uni_com}
There exists a complex of $S_\bullet$-modules $\EE = \EE(\lambda)$, such that for any $(m, n)$ the first layer of the spectral sequence associated to the Lie subalgebra
$\gl(m, n - 1) \oplus \gl(1) \into \gl(m, n)$ and coefficients $\Sigma^\lambda(V) \tensor \Sigma^\lambda(V^*)$ has the form
$$
E_1^{pq} \ \isom\ \bigoplus_{\mu\atop i + j = q} H^i(\gl(m, n - 1), \Sigma^\mu(W) \tensor \Sigma^\mu(W^*)) \tensor H^j(\gl(1), \k) \tensor \EE^p_\mu.
$$
where $W$ is the standard representation of $\gl(m, n - 1)$, and the sum is over $\mu \in \cH_{m,n - 1}$.

Moreover, the differentials in $E_1$ are induced by the differentials in $\EE$, so the second layer $E_2$ (without differentials) has the form
$$
E_2^{pq} \ \isom\ \bigoplus_{\mu\atop i + j = q} H^i(\gl(m, n - 1), \Sigma^\mu(W) \tensor \Sigma^\mu(W^*)) \tensor H^j(\gl(1), \k) \tensor H^p(\EE_\mu).
$$
\end{lemma}
\qed

Because of the universal nature of the complex $\EE$ we will be able to obtain information about $H^\bullet(\EE)$ by comparing spectral sequences for different values
of $(m, n)$. First, we will need the following simple lemma.

\begin{lemma}
\label{lemma_sp_seq}
Consider a spectral sequence $E_\bullet$ concentrated in the quadrant with $p, q \ge 0$, such that $E_2^{0,q}$ is an algebra isomorphic to $A \tensor \k[e_1]$ for some algebra $A$,
and $\deg e_1 = 1$.
Assume that the second layer is a free $A \tensor \k[e_1]$-module, with generators in degrees $(p, 0)$, and the differential on $E_2$ is compatible with the module structure.
\begin{enumerate}[label=\alph*)]
\item If $E_\infty^{0,q} \isom E_2^{0,q}$ and $E_\infty^{pq} = 0$ for $p > 0$, then $E_2^{pq} = 0$ for $p > 0$.
\item If $E_\infty^{0,q} \isom A \into E_2^{0,q}$ and $E_\infty^{pq} = 0$ for $p > 0$, then
$$
E_2^{pq} \isom A \tensor \k[e_1] \tensor \k[c_1]
$$
as a $A \tensor \k[e_1]$-module. Here $\deg c_1 = (2, 0)$ and the differential in $E_2$ sends $e_1$ to $c_1$.
\end{enumerate}
\end{lemma}

\proof
For part (a) let us assume that $E_2^{pq} \neq 0$ for some $p > 0$, and let us denote $p_0$ the minimal such $p$. By our assumption $E_2$ is generated by elements in degrees $(p, 0)$,
therefore we must have $E_2^{p_0, 0} \neq 0$. Since these elements do not survive to $E_\infty$ and $p_0$ is minimal, we must have non-zero differentials starting from the first column. But this is
impossible since $E_\infty^{0,q} \isom E_2^{0,q}$.

For part (b) observe that since $e_1$ doesn't survive until $E_\infty$ it must by killed by some differential, however, since our spectral sequence has non-zero terms only for $p, q \ge 0$,
we see that the differential of $E_2$ doesn't vanish on $e_1$. Denote by $c_1 \in E_2^{2,0}$ its image. Since by assumption the differential is compatible with the module structure
and $A$ is contained in the kernel of the differential, this completely determines restriction of $d_2$ to the first column. Furthermore, since
$E_2$ is a free $A \tensor \k[e_1]$-module, we see that the image $\Im d_2\from E_2^{0\bullet} \to E_2^{2\bullet}$ is identified with $Ac_1 \subset (A \tensor \k[e_1]) c_1$, and
it is surjective on $E_2^{2,0}$.

The element $e_1 c_1$ does not survive to $E_\infty$, therefore as before, $d_2$ doesn't vanish on it and we identify $d_2(e_1 c_1)$ with $c_1^2$. Repeating this argument we obtain the
required isomorphism.

\qed

These two lemmas allow us to establish the following result concerning the universal complex $\EE(\lambda)$.

\begin{lemma}
\label{lemma_EE}
Let $\EE = \EE(\lambda)$, then
\begin{enumerate}[label=\alph*)]
\item for all $p > 0$ we have $H^{2p}(\EE) \isom L_\mu$, for some $\mu$ with $\rmht(\mu) = \rmht(\lambda) + 1$,
\item $H^0(\EE) \isom L_{\wbar\lambda}$, where $\wbar\lambda$ is the truncation of $\lambda$,
\item for all $p \ge 0$ we have $H^{2p}(\EE) \isom L_{(\lambda / p + 1)}$, where $(\lambda / k)$ denotes the diagram obtained from $\lambda$ by adding one box
in the first $k - 1$ columns and removing $k$'th column. In other words
$$
(\lambda / k)' = (\lambda'_1 + 1, \ldots, \lambda'_{k - 1} + 1, \lambda'_{k + 1}, \ldots).
$$
\end{enumerate}
\end{lemma}

\proof
Consider Lie superalgebra $\g = \gl(k, 1)$ and its subalgebra $\h = \gl(k) \oplus \gl(1)$. Now, as usual, let $W$ be the standard representation of $\gl(k)$
then for any diagram $\mu$ with $\rmht(\mu) \le k$ we have
$$
H^\bullet(\gl(k), \Sigma^\mu(W) \tensor \Sigma^\mu(W^*)) \isom \k[e_1, \ldots, e_{2k - 1}].
$$
We also have $H^\bullet(\gl(1)) \isom \k[e_1]$. Using lemma \ref{lemma_uni_com} we see that the spectral sequence for the Lie subalgebra $\h \into \g$ satisfies
conditions of lemma \ref{lemma_sp_seq}.

First, let $k \ge \rmht(\lambda) + 1$. As we already saw, all the diagrams $\mu$ appearing in the universal complex $\EE$ have $\rmht(\mu) \le \rmht(\lambda) + 1$.
Using theorem \ref{thm_cohom_glm1} we find that the spectral sequence converges to $\k[e_1, \ldots, e_{2k-1}]$. Therefore, we are in the situation of lemma
\ref{lemma_sp_seq}(b).

On the other hand, let $k = \rmht(\lambda)$. Then again using theorem \ref{thm_cohom_glm1} we find that the spectral sequence converges to $\k[e_1, \ldots, e_{2k - 1}, e'_1]$,
so we are in the situation of lemma \ref{lemma_sp_seq}(a). The difference between these two spectral sequences is that in the latter case we have a restriction imposed on the
diagrams $\mu$ by the dimension of $W$. Namely, we lose diagrams $\mu$ with $\rmht(\mu) > k$. Therefore, diagrams $\mu$ contributing to $H^p(\EE)$ for $p > 0$ must
be of height $\rmht(\lambda) + 1$. Combining this with isomorphism from lemma \ref{lemma_sp_seq}(b) we prove part (a).

For part (b), notice that since $H^0(\gl(k, 1), \Sigma^\lambda(V) \tensor \Sigma^\lambda(V^*)) = \k$ for any $k$, the cohomology $H^0(\EE)$ is isomorphic to a simple module.
Furthermore, from the construction of the first layer $E_1$ we see that it must by isomorphic to $L_\mu$ for the smallest possible diagram $\mu$ in the decomposition (\ref{equ_VEW_coeff}).
It is straightforward to see that this diagram is precisely the truncated diagram $\wbar\lambda$.

Clearly, $(\lambda / 1) = \wbar\lambda$, so part (c) for $p = 0$ is just a reformulation of part (b).
Let us consider the rest of the diagrams $\mu$ contributing to $H^\bullet(\EE(\lambda))$. According to part (a) they are of height $\rmht(\lambda) + 1$, so they have $\rmht(\lambda) + 1$
boxes in the first column. This is only possible if in the process of obtaining $\mu$ from $\lambda$ no box was removed from the first column and exactly one box was added to it.
This implies that for any diagram $\mu$ of height $\rmht(\lambda) + 1$ appearing in the universal complex $\EE(\lambda)$ we have the truncated diagram $\wbar\mu$ appearing in $\EE(\wbar\lambda)$.
More precisely, for any such $\mu$ and $p \ge 2$ we have $\EE(\lambda)^p_\mu \isom \EE(\wbar\lambda)^{p-2}_{\wbar\mu}$ and therefore
$$
H^p(\EE(\lambda)_\mu) \isom H^{p-2}(\EE(\wbar\lambda)_{\wbar\mu}).
$$
Thus, we reduced the question to the structure of the universal complex $\EE(\wbar\lambda)$, for a diagram $\wbar\lambda$ that contains one less column than $\lambda$.
Using induction on the number of columns we immediately see that this identification gives us the required isomorphism of part (c). The base of induction is the diagram
$\lambda$ of size zero, in which case the statement follows from the decomposition (\ref{equ_VEW_ext}).

\qed

We are now ready to prove the main theorem of this section.

\begin{theorem}
Let $\g = \gl(m, n)$ with $m \ge n \ge 0$, and
$$
\lambda \in \cH_{m - n + k, k} - \cH_{m - n + k - 1, k - 1}
$$
for some $0 \le k \le n$. Then
$$
H^\bullet(\g, \Sigma^\lambda(V) \tensor \Sigma^\lambda(V^*)) \ \isom\ \k[e_1, \ldots, e_{2m - 1}] \tensor \k[e'_{2(n - k) + 1}, \ldots, e'_{2n - 1}],
$$
where $\deg e_i = \deg e'_i = i$. The generators $e_1,\ldots, e_{2m-1}$ are the images of the standard generators of $H^\bullet(\gl(m), \k)$ under the composition
\begin{equation}
\label{equ_classes_e}
\begin{tikzcd}[column sep=4em]
H^\bullet(\gl(m), \k) & H^\bullet(\g, \k) \ar[l, "\isom"', "\mathrm{res}"] \ar[hook,r,"\mathrm{coev}"'] & H^\bullet(\g, \Sigma^\lambda(V) \tensor \Sigma^\lambda(V^*)).
\end{tikzcd}
\end{equation}
Here $\mathrm{res}$ is the map induced by restriction to the Lie subalgebra $\gl(m) \into \g$, and $\mathrm{coev}$ is induced by the coevaluation map
$\k \into \Sigma^\lambda(V) \tensor \Sigma^\lambda(V^*)$.
\end{theorem}

\proof
We consider Lie subalgebra $\h  = \gl(m, n - 1) \oplus \gl(1) \into \gl(m, n)$ and prove the theorem by induction on $n$. When $n = 0$ this is the classical purely even case
and the result is well known. When $n = 1$ the statement of the theorem is a reformulation of theorem \ref{thm_cohom_glm1}. From now on, we will assume that $n \ge 2$ and the theorem holds for all $\gl(m, n')$ with $n' < n$. According to
lemma \ref{lemma_uni_com} we have the spectral sequence
$$
E_2 \ \isom\ \bigoplus_{\mu} H^\bullet(\gl(m, n - 1), \Sigma^\mu(W) \tensor \Sigma^\mu(W^*)) \tensor \k[e''_1] \tensor H^\bullet(\EE_\mu) \Rightarrow H^\bullet(\g, \Sigma^\lambda(V) \tensor \Sigma^\lambda(V^*)).
$$

\begin{nparagraph}[Case $k = 0$.]
First, we consider the case when $k = 0$ separately. The condition $\lambda \in \cH_{m - n, 0}$ is equivalent to $\rmht(\lambda) \le m - n$.  From lemma \ref{lemma_EE} we see that all diagrams
$\mu$ contributing to the second layer of the spectral sequence have $\rmht(\mu) \le m - n + 1$, i.e. $\mu \in \cH_{m - (n - 1), 0}$, hence by the inductive assumption
$$
H^\bullet(\gl(m, n - 1), \Sigma^\mu(W) \tensor \Sigma^\mu(W^*)) \ \isom\ \k[e_1, \ldots, e_{2m - 1}].
$$
Furthermore, again by lemma \ref{lemma_EE} we have
$$
E_2 \ \isom\ \k[e_1, \ldots, e_{2m - 1}] \tensor \k[e''_1] \tensor \k[c_1]
$$
and the differential on $E_2$ sends $e'_1$ to $c_1$. Therefore, we find that the spectral sequence converges to
$$
H^\bullet(\g, \Sigma^\lambda(V) \tensor \Sigma^\lambda(V^*)) \ \isom\ \k[e_1, \ldots, e_{2m - 1}].
$$

The statement regarding classes $e_i$ is a tautology for $n = 0$. Assume that it hold for all $n' < n$. The coevaluation map $\k \to \Sigma^\lambda(V) \tensor \Sigma^\lambda(V^*)$ induces
a morphism of spectral sequences. By inductive assumption this morphism is an isomorphism on the second layer, therefore it also induces isomorphism
$H^\bullet(\g, \k) \isom H^\bullet(\g, \Sigma^\lambda(V) \tensor \Sigma^\lambda(V^*))$.

This completes the proof in the case $k = 0$.
\end{nparagraph}

\begin{nparagraph}[Case $k > 0$.]
According to lemma \ref{lemma_EE}(c) the contribution to the second layer of the spectral sequence in the column $2p$ comes from the diagram $(\lambda / p + 1)$.
To simplify notation, let us write
$$
\cH^\circ(n, k) = \cH_{m - n + k, k} - \cH_{m - n + k - 1, k - 1}.
$$
We omit $m$ in the notation since within the scope of this proof the number of even variables $m$ never changes. It is straightforward to check that since $\lambda \in \cH^\circ(n, k)$
for $0 \le p \le k - 1$ we have
$$
(\lambda / p + 1) \in \cH^\circ(n - 1, k - 1),
$$
and for $p \ge k$
$$
(\lambda / p + 1) \in \cH^\circ(n - 1, k).
$$

Therefore, by inductive assumption for $0 \le p \le k - 1$ the $2p$'th column of $E_2$ is isomorphic to
$$
E_2^{2p,\bullet} \isom \k[e_1, \ldots, e_{2m-1}] \tensor \k[e'_{2(n-k)+1}, \ldots, e'_{2n-3}] \tensor \k[e''_1].
$$

For $p \ge k$ we consider two cases. First, if $k < n$, then again by inductive assumption we have
$$
E_2^{2p,\bullet} \isom \k[e_1, \ldots, e_{2m-1}] \tensor \k[e'_{2(n-k)-1}, \ldots, e'_{2n-3}] \tensor \k[e''_1].
$$
If on the other hand, $k = n$, then $(\lambda / p + 1) \in \cH^\circ(n - 1, n)$, hence the Schur functor $\Sigma^{(\lambda / p + 1)}(W) = 0$, and all
columns in $E_2$ starting from column $2k$ vanish.

Since the differential on $E_2$ sends generator $e''_1 \in E_2^{2p,1}$ to the basis element in $E_2^{2(p+1),0}$ we see that on the third layer the spectral sequence
has only two non-zero columns: for $p = 0$ and either for $p = 2k$ if $k < n$ or for $p = 2k - 2$ if $k = n$. Specifically,
$$
E_3^{0,\bullet} \isom \k[e_1, \ldots, e_{2m-1}] \tensor \k[e'_{2(n-k)+1}, \ldots, e'_{2n-3}],
$$
and if $0 < k < n$, then
$$
E_3^{2k,\bullet} \isom \k[e_1, \ldots, e_{2m-1}] \tensor \k[e'_{2(n-k)+1}, \ldots, e'_{2n-3}] \tensor \k e'_{2(n-k)-1},
$$
and finally if $k = n$, then
$$
E_3^{2(k-1),\bullet} \isom \k[e_1, \ldots, e_{2m-1}] \tensor \k[e'_{1}, \ldots, e'_{2n-3}] \tensor \k e''_1.
$$

Let us show that starting from $E_3$ all differentials in the spectral sequence vanish. First, consider generators $e_i$. The coevaluation map $\k \to \Sigma^\lambda(V) \tensor \Sigma^\lambda(V^*)$
induces a map from the spectral sequence for the trivial coefficients $F_\bullet$ to our spectral sequence $E_\bullet$. By inductive assumption the classes $e_i$ in the first column of $F_\bullet$
map to corresponding classes $e_i$ in $E_\bullet$. As we have seen for the trivial coefficients all the differentials vanish on $e_i$, hence they must also vanish in $E_\bullet$.

Now consider generators $e'_j$, and assume first that $k < n$. The differential can only be non-zero on the layer $E_{2k}$ and send generator $e'_j$ to $E_{2k}^{2k,j - 2k + 1}$. However,
since $j \le 2n-3$ we have $j - 2k + 1 \le 2(n - k) - 2$ and column $2k$ has non-zero terms only for $q \ge 2(n - k) - 1$.

Finally, for $k = n$, the differential can only be non-zero on the layer $E_{2(n-1)}$ and send generator $e'_j$ to $E_{2(n-1)}^{2(n-1),j - 2n + 3}$. Again, since $j \le 2n-3$ we have
$j - 2n + 3 \le 0$ but all non-zero terms are in degree $q \ge 1$.

Observe, that for $k < n$ the total degree of the generator $e'_{2(n-k)-1}$ in the column $2k$ is $(2n - 1)$, and similarly for $k = n$ the total degree of $e''_1$ in the column
$2(n - 1)$ is again $(2n - 1)$. By renaming this generator $e'_{2n-1}$ we obtain the required isomorphism. The identification (\ref{equ_classes_e}) immediately follows from the
previous discussion.

This concludes the proof of the theorem.

\qed

\end{nparagraph}

\vskip 5em
\section{Cohomology of $\V_{m,n}$}
\label{sec_Vmn}
The calculation of cohomology $H^\bullet(\V_{m,n}, \k)$ follows the general argument originally developed for the classical case of $\V_{m,0}$ by Gelfand and Fuchs
with some refinements that were needed to apply it to $\V_{m,1}$. Here we briefly recall the major steps of this procedure, for details we refer to \cite{Fuks} and \cite{gfn1}.

\begin{nparagraph}
Consider Lie subalgebra $\gl(m, n) \into \V_{m,n}$, and let $V$ be the standard representation of $\gl(m, n)$. The continuous dual space
$$
\Hom(\V_{m,n}, \k) \ \isom\ \bigoplus_{i \ge 0} \left( S^i(V^*) \tensor V \right).
$$
Therefore, in the spectral sequence for the Lie subalgebra $\gl(m, n) \into \V_{m,n}$
$$
E_1^{pq} = H^q(\gl(m,n), \Hom(\Lambda^p(\V_{m,n} / \gl(m,n)), \k)) \Rightarrow H^{p+q}(\V_{m,n}, \k),
$$
the coefficients of the cohomology groups of $\gl(m,n)$ can be written as
$$
\bigoplus_{\sum p_i = p} \Lambda^{p_i} \left(S^i(V^*) \tensor V\right),
$$
where $i \ge 0$ and $i \neq 1$. This can be simplified by observing that the contributions to the first layer of the spectral sequence can only come
from terms of the form
$$
\Lambda^p(V) \tensor \Lambda^p(S^2(V^*) \tensor V).
$$
By expanding the second exterior power and using calculus of Schur functors one then shows that the first layer of the spectral sequence has the form
\begin{equation}
\label{equ_gf_E1}
E_1^{2p,q} \ \isom\ \bigoplus_{|\lambda| = p} H^\bullet(\gl(m,n), \Sigma^{\wtilde\lambda}(V) \tensor \Sigma^{\wtilde\lambda}(V^*)),
\end{equation}
where $\wtilde\lambda$ is obtained from $\lambda$ by adding to it one more column with $|\lambda|$ boxes in it, in other words $\wtilde\lambda_i = \lambda_i + 1$ for
$1 \le i \le p$.

\end{nparagraph}

\begin{theorem}
For any $m \ge n \ge 0$ we have an isomorphism
$$
H^\bullet(\V_{m,n}, \k) \ \isom\ H^\bullet(\SS^{2n} X_{2(m - n)}, \k).
$$
\end{theorem}

\proof
First of all notice that if $\lambda$ is a diagram in (\ref{equ_gf_E1}) then its transposed $\lambda'$ is a diagram appearing in the similar spectral sequence
for the Lie subalgebra $\gl(n - 1, m + 1) \into \V_{n - 1, m + 1}$. To simplify notation we put
$$
\cH^\circ(m, n, k) = \cH_{m - n + k, k} - \cH_{m - n + k - 1, k - 1}.
$$
Clearly, if $\lambda \in \cH^\circ(m, n - 1, k)$ for $k \ge 1$, then $\wtilde\lambda \in \cH^\circ(m, n, k + 1)$. If $k = 0$, then there are two possibilities:
if $|\lambda| \le m - n$, then $\wtilde\lambda \in \cH(m, n, 0)$, otherwise $\wtilde\lambda \in \cH^\circ(m, n, 1)$.

Denote by $\hat\lambda$ the diagram obtained from $\lambda$ by adding one more row with $|\lambda|$ boxes, in other words
$$
\hat\lambda = (|\lambda|, \lambda_1, \lambda_2, \ldots) = \left(\wtilde{\lambda'}\right)'.
$$
If $\lambda \in \cH^\circ(m, n - 1, k)$ for any $k \ge 0$, then $\hat\lambda \in \cH^\circ(m + 1, n - 1, k)$.

\begin{nparagraph}
Let us denote by $F_\bullet$ the spectral sequence (\ref{equ_gf_E1}) for $\V_{n - 1, m + 1}$. On the first layer $F_1$ the term corresponding to the diagram
$\lambda \in \cH^\circ(m, n - 1, k)$ for $0 \le k \le n - 1$ is isomorphic to
$$
(F_1)_\lambda \isom \k[e_1, \ldots e_{2m+1}] \tensor \k[e'_{2(n-k)-1}, \ldots e'_{2n-3}].
$$
And the corresponding term in the spectral sequence $E_1$ is isomorphic to
$$
(E_1)_\lambda \isom \begin{cases}
\k[e_1, \ldots e_{2m-1}] \tensor \k[e'_{2(n-k)-1}, \ldots e'_{2n-1}],&\text{if $|\lambda| > m - n$},\\
\k[e_1, \ldots e_{2m-1}],&\text{if $|\lambda| \le m - n$}.
\end{cases}
$$

Since $n - 1 < m + 1$ the cohomology of $\V_{n-1,m+1}$ is covered by theorem \ref{thm_known_gf}(b). Therefore, the spectral sequence $F_\bullet$ converges to
$H^\bullet(S^{2m + 1}, \k) = \k[e_{2m+1}]$. Moreover, as was shown in \cite{Fuks-super} this class $e_{2m+1}$ maps to the corresponding class in $H^\bullet(\V_{0,m+1}, \k)$
under the restriction map to the Lie subalgebra $\V_{0,m+1} \into \V_{n-1,m+1}$. And from the discussion in \cite{gfn1} section 3, it follows that this class further maps
to $e_{2m+1} \in H^\bullet(\gl(n - 1, m + 1), \k)$ under the restriction map to the Lie subalgebra $\gl(n - 1, m + 1) \into \V_{n-1,m+1}$. Hence, all the differentials
in $F_\bullet$ vanish on the generator $e_{2m+1}$ and from the degree considerations $e_{2m+1}$ doesn't appear in the image of the differentials of any other generator
$e_i$ or $e'_j$. So we have the sub-spectral sequence of $F_\bullet$ that we will denote by $G_\bullet$, such that
$$
F_\bullet \isom G_\bullet \tensor \k[e_{2m+1}].
$$
The component of $G_1$ corresponding to a diagram $\lambda \in \cH^\circ(m, n - 1, k)$ for $0 \le k \le n - 1$ is
$$
(G_1)_\lambda \isom \k[e_1, \ldots e_{2m-1}] \tensor \k[e'_{2(n-k)-1}, \ldots e'_{2n-3}],
$$
and $G_\bullet$ converges to $\k$ (in degree $0$).

\end{nparagraph}

\begin{nparagraph}
Let us compare spectral sequences $E_\bullet$ and $G_\bullet$. We introduce an intermediate spectral sequence $\wtilde E_\bullet$ by adding to $E_\bullet$ the ``missing'' classes
$e_{2n - 1}$ to all the small diagrams $\lambda$ with $|\lambda| \le m - n$. The differentials in the spectral sequence $E_\bullet$ can be described as follows. For generators
$e_{2i-1}$ the differentials $d_r$ vanish up to layer $r = 2i$, and on the layer $E_{2i}$ they send
$$
(e_{2i-1} z_\lambda) \mapsto \sum_{\mu \in \lambda \cdot c_i} z_\mu,
$$
where sum is taken over diagrams $\mu$
in the decomposition of the product of $\lambda$ and the Chern class $c_i$, i.e. $\mu$ is obtained from $\lambda$ by adding $i$ boxes such that no more than one box added in each row. Here $z_\lambda$
denotes the generator of the component $(E_1)_\lambda$. Similarly, for generators $e'_{2j-1}$ differentials vanish up to layer $r = 2j$ and on that layer they send
$$
e'_{2j-1} z_\lambda \mapsto \sum_{\nu \in \lambda \cdot s_j} z_\nu,
$$
where $\nu$ is a diagram in the decomposition of the product of $\lambda$ and Segre class $s_j$, i.e. $\nu$ is obtained from $\lambda$ by adding $j$ boxes, such that
no more than one box is added in each column.

We define spectral sequence $\wtilde E_\bullet$ by setting
$$
(\wtilde E_1)_\lambda = \begin{cases}
(E_1)_\lambda,&\text{if $|\lambda| > m - n$},\\
(E_1)_\lambda \tensor \k[e_{2n-1}],&\text{if $|\lambda| \le m - n$}.
\end{cases}
$$
The differentials are defined as described above.

In fact one can construct a filtered complex for this spectral sequence $\wtilde E_\bullet$. We start from the cochain complex $C^\bullet = C^\bullet(\V_{m,n}, \k)$ with the filtration
induced by the Lie subalgebra $\gl(m,n)$ as described in \ref{par_sp_seq}. In every degree $p \ge 0$ this is a bounded filtration of $C^p$, therefore
the filtration of the entire complex is both complete and cocomplete. In such case we can construct a bicomplex $B^{\bullet\bullet}$, so that its totalization equipped
with one of the natural filtrations of the bicomplex (say in the vertical direction) is filtered quasi-isomorphic to the cochain complex $C^\bullet$.
This can be seen as a special case of Koszul duality between filtered complexes, that are identified via Rees construction with complexes of (flat) $\k[u]$-modules, and
on the dual side complexes of $\k[\epsilon]$-modules, that are identified with bicomplexes, where horizontal differential is given by $d$ and vertical differential
by the action of $\epsilon$ (for a brief summary we refer to \cite{bbta} section 1.2, for a detailed discussion see for example \cite{Pos}).

Now, since the totalization $\Tot B^{\bullet\bullet}$ is filtered quasi-isomorphic to $C^\bullet$ the first layers of the corresponding spectral sequences are isomorphic. We construct
bicomplex $\wtilde B^{\bullet\bullet}$ by putting
$$
\wtilde B^{p\bullet} \isom \begin{cases}
B^{p\bullet},&\text{if $p > m - n$},\\
B^{p\bullet} \tensor \k[e_{2n-1}],&\text{if $p \le m - n$}.
\end{cases}
$$
We put both horizontal and vertical differentials to be zero on $e_{2n-1}$. We will denote $\wtilde C^\bullet$ the totalization $\Tot \wtilde B^{\bullet\bullet}$ equipped with the
filtration in the vertical direction.

Notice that $\wtilde E_1 \isom G_1 \tensor \k[e_{2n-1}]$. In the spectral sequence $G_\bullet$ the roles of generator classes $e$ and $e'$ are reversed, however, the diagrams
appearing in $G_\bullet$ are the transposes of those appearing in $E_\bullet$, therefore the differentials in $G_\bullet$ have the same description as above.
Therefore, since the spectral sequence $G_\bullet$ converges to $\k$ all potential targets for differential starting from
the new class $e_{2n-1}$ are already killed in $G_\bullet$. Hence, we find that $\wtilde E_\bullet$ converges to $\k[e_{2n-1}]$.

\end{nparagraph}

\begin{nparagraph}
Finally, consider the short exact sequence $C^\bullet \to \wtilde C^\bullet \to Q^\bullet \tensor \k e_{2n - 1}$. Here the spectral sequence for $Q^\bullet$ has only contributions
from diagrams $\lambda$ with $|\lambda| \le m - n$ and
$$
(Q_1)_\lambda \isom \k[e_1, \ldots e_{2m-1}].
$$
This, as in the classical case, is isomorphic to the spectral sequence for the fiber product
$$
X_{2(m-n)} = \sk_{2(m-n)} BGL(m) \times_{BGL(m)} EGL(m).
$$

Since $\wtilde E_\bullet$ converges to $\k[e_{2n-1}]$, from the long exact sequence we find that
\begin{align*}
H^0(\V_{m,n}, \k) &= \k,\\
H^i(\V_{m,n}, \k) &= 0,\quad \text{for $1 \le i \le 2n$}\\
H^i(\V_{m,n}, \k) &\isom H^{i-2n}(X_{2(m-n)},\k), \quad \text{for $i > 2n$}.
\end{align*}

\qed

\end{nparagraph}

\vfill\eject

\end{document}